\documentclass[10pt]{amsart}
\usepackage{amsfonts,amssymb,latexsym}
\usepackage{pb-diagram}
\usepackage{amscd,amsmath}

\usepackage[arrow,matrix,curve]{xy}\SilentMatrices
\def\xyma{\xymatrix@M.7em}

\newtheorem{theorem}{Theorem}[section]
\newtheorem{lemma}[theorem]{Lemma}
\newtheorem{definition}[theorem]{Definition}
\newtheorem{example}[theorem]{Example}

\def\g{\gamma}
\def\G{\Gamma}

\def\S{\Sigma}
\def\Z{\mathbb{Z}}
\def\R{\mathbb{R}}
\def\Q{\mathbb{Q}}

\begin{document}
\title{A new look at finitely generated metabelian groups}
\author{Gilbert Baumslag}
\address{Dept. of Computer Science,
City College of New York, Convent Avenue and 138th Street, New York, NY 10031}
\email{gilbert.baumslag@gmail.com}
\thanks{The research of the first author is supported by Grant CNS 111765.}

\author{Roman Mikhailov}
\address{Steklov Mathematical Institute, Gubkina 8, 119991 Moscow, Russia and Institute for Advanced Study, Princeton, NJ, USA}
\email{romanvm@mi.ras.ru}
\urladdr{http://www.mi.ras.ru/\~{}romanvm/pub.html}
\thanks{The research of the second author is
supported by the National Science Foundation under agreement No.
DMS-0635607. Any opinions, findings and conclusions or
recommendations expressed in this material are those of the
authors and do not necessarily reflect the views of the National
Science Foundation.}

\author{Kent E. Orr}
\address{Dept of Mathematics, Indiana University, Bloomington IN 47405}
\email{korr@indiana.edu}
\thanks{The third author thanks the National Science Foundation, Grant 707078, and the Simons Foundation, Grant 209082, for their support.}

\date{\today}

\maketitle
\begin{abstract}
A group is metabelian if its commutator subgroup is abelian. For  finitely generated metabelian groups, classical commutative algebra, algebraic geometry and
 geometric group theory, especially the latter two subjects, can be brought to bear on their study. The object
 of this paper is to describe some of the new ideas and open problems that arise.
\end{abstract}

\section{Introductory remarks}

There was a conference at the City College of New
York on March 17 and 18 of 2011, funded by the National Science
Foundation, entitled {\em Finitely Presented Solvable Groups.} The
first named author of this paper discussed how geometric group
theory and algebraic geometry might play a role in better
understanding finitely generated metabelian groups, in particular
the isomorphism problem. In the past 18 months the present authors
have introduced some new ideas which more closely tie algebraic
geometry to the isomorphism problem, in particular, and to the
study of these finitely generated metabelian
groups~\cite{Baumslag-Mikhailov-Orr:2012-1,Baumslag-Mikhailov-Orr:2012-2}. A
key idea is that of {\em para-equivalence of finitely generated
metabelian groups}, in which localization of modules plays a
fundamental and simplifying role, as it does in algebraic
geometry. The object of the present note is to provide a quick
overview of the general subject and to explain how these ideas
come into play. Subsequent papers will describe this in more
detail.

Geometric group theory in the form of hyperbolic groups was
introduced by Gromov in 1983 \cite{Gromov:1983-1}. Gromov's idea was
to view infinite groups as geometric objects, and with his
extraordinary power and insight he revolutionized part of
combinatorial group theory. This founded the subject of {\em
geometric group theory.} In fact, many of the ideas go back to
Dehn  in 1910, 1911, Tartakovski~\cite{Tartakovski:1949-1},~\cite{Tartakovski:1949-2},
\cite{Tartakovski:1949-3}, Lyndon, Weinbaum,
Greendlinger, Lipschutz, Schupp (see the book~\cite{Lyndon-Schupp:1977-1}
for these and a host of references), and Rips~\cite{Rips:1982-1}.
These hyperbolic groups are a far cry from solvable groups, and
more closely resemble free groups. In fact solvable groups have
been largely excluded from geometric group theory.  Recently, some
attention has been paid to metabelian groups - witness the work of
Farb and Mosher \cite{Farb-Mosher:2000-1}, Eskin and Fisher and Whyte \cite{Eskin-Fisher-Whyte:2007-1}, \cite{Eskin-Fisher-Whyte:2007-2}.

The hyperbolic groups of Gromov are finitely presented, unlike
many finitely generated metabelian groups. However, as we will
point out in due course, when viewed as objects in the category of
metabelian groups, they are finitely presented and geometric ideas
might, if properly developed, play a part in their study.
Geometric ideas of a rather different kind do arise naturally in
the case of metabelian groups through commutative algebra. This
suggest that major progress may arise through ideas borrowed from
algebraic geometry.

\section{Philip Hall's approach}

In 1954, Philip Hall began a systematic study of
finitely generated metabelian groups~\cite{Hall:1954-1} . We
recall his approach here. To this end, we shall use the following
notation. Given a group $G$, we denote the conjugate $g^{-1}ag$ of
an element $a\in G$ by the element $g$ by $a^g$, the commutator
$g^{-1}a^{-1}ga$ of $g$ and $a$ by $[g,a]$ and the subgroup of a
group $G$ generated by the commutators $[x,y]=x^{-1}y^{-1}xy$,
where $x$ and $y$ range over the elements of $G$,  by $[G,G]$.
Then $G$ is termed  metabelian if $[G,G]$ is abelian. Thus a
metabelian group $G$ is an extension of an abelian normal subgroup
$A$ by an abelian group $G/A=Q$.

Hall first observed that one can view this subgroup $A$ as a $\Z[Q]$-module.  Precisely, $Q$ acts on $A$ on the right by conjugation:
\[
a(gA)=g^{-1}ag\ (a\in A, g\in G).
\]
So the integral group ring $R=\Z[Q]$  of $Q$ acts on $A$ which, writing it additively,
then becomes a right $R$-module:
\[
a(\Sigma_{i=1}^nc_ig_iQ)=\Sigma_{i=1}^nc_i(g_i^{-1}ag_i)\, (g_i \in G).
\]
If $G$  is finitely generated, then $Q$ is a finitely generated
abelian group. It  then turns out that $A$ is a finitely generated
$R$-module and so one can avail oneself of the structure theory of
classical commutative algebra. In particular, $A$ is Noetherian.

\section{A few properties of finitely generated metabelian groups}

Philip Hall's approach to the study of finitely
generated metabelian groups described above in $\mathsection 2$
gave rise to an entire body of results. We describe
some of these below.  Items 1, 2, 3 and 6 are due to Hall. A
general reference to much of what is described in this section can
be found in the book \cite{Lennox-Robinson:2004-1}.

\begin{enumerate}
\item  Finitely generated metabelian groups satisfy the maximal condition for normal
subgroups, that is, any properly ascending chain of normal subgroups is finite. This follows
from the fact, adopting the notation above,  that $A$ is Noetherian.
\item Every finitely generated metabelian group $G$ is finitely presentable in the category
of metabelian groups, i.e.,  $G$ can be defined in terms of finitely many generators and
finitely many relations together with all relations of the form $[[x,y],[z,w]]=1$. We write a group presentation in the category of metabelian groups using double angle brackets, as follows:
\[
G=\langle\langle x_1,\dots,x_m\ |\
r_1=1,\dots,r_n=1\rangle\rangle.
\]
\item  It follows from (2) that there are only countably many isomorphism classes of finitely generated
metabelian groups.
\item Algorithmically, finitely generated metabelian groups are well behaved. So for example, the word and conjugacy
problem are solvable.
\item There is an algorithm to decide if a finitely generated metabelian group is residually nilpotent,
i.e., if the intersection of the terms of the lower central series is trivial.
\item Finitely generated metabelian groups are residually finite, i.e., the normal subgroups of
finite index have trivial intersection.
\item Although finitely generated metabelian groups are not necessarily finitely presented,
every finitely presented metabelian group can be embedded in a finitely presented metabelian group.
\item Bieri-Strebel have introduced a geometric invariant which distinguishes the finitely presented metabelian groups
from the others~\cite{Bieri-Strebel:1981-1}. It is commonly believed that the computation of this invariant is algorithmic.
\item Groves-Miller~\cite{Groves-Miller-Charles:1986-1} have solved the isomorphism problem for free metabelian groups, i.e., for metabelian groups
which can be presented in the form
\[
G=\langle\langle x_1,\dots,x_m\rangle\rangle.
\]
\item As of now the isomorphism problem for finitely generated metabelian groups remains, perhaps,
the most intriguing open problem about these groups. At first sight it seems to be connected to Hilbert's
Tenth problem which suggests that this problem is algorithmically undecidable. We shall not touch on
this aspect of the isomorphism problem but go in a different direction. We will say more about the isomorphism
problem in the course
of further discussion.
\end{enumerate}

\section{The Bieri-Strebel invariant}

Let $Q$ be a finitely generated abelian group. As
previously noted, for any extension $G$ of a (not necessarily
finitely generated) abelian group $A$ by $Q$, the conjugation
action of $G$ on $Q$ makes $A$ into a $\Z[Q]$--module. It is easy to
see that $G$ is finitely generated if and only if $A$ is finitely
generated as a $\Z[Q]$--module. A homomorphism of $Q$ to the additive
group $\R$ of real numbers is called a {\em valuation} of $Q$.
Associated to each such valuation $v$ one has the submonoid of $Q$
\[
Q_v = \{ q\in Q\mid v(q)\ge 0\}.
\]
For valuations $v,v'$ we write $v\sim v'$ if and only if there exists
$\lambda > 0$ such that $v(q)=\lambda v'(q)$ for all $q\in Q$.
Let $n$ be the torsion-free rank of $Q$.
Then ${Hom} (Q,\mathbb R) \cong \mathbb R^n$, and
there is an obvious identification between the set of equivalence classes
of nontrivial valuations of $Q$ and the $(n-1)$-sphere $S^{n-1}$.

Let $A$ be a finitely generated $\Z[Q]$--module.
We can view $A$ as a module over the commutative ring
$\mathbb Z[Q_v] \leq \mathbb Z [Q]$.
Define $\Sigma_A$ to be the set of $\sim$ classes of valuations $v$ on $Q$
such that $A$ is finitely generated as a $\mathbb Z [Q_v] $--module.

The module $A$ is said to be {\em tame} if $\S_A \cup -\S_A = S^{n-1}$,
in other words, for every valuation $v$ of $Q$, either $A$
is finitely generated as a $\Z[Q_v]$--module, or else it is finitely
generated as a $\Z[Q_{-v}]$--module.

Then Bieri and Strebel have proved the following theorems~\cite{Bieri-Strebel:1981-1}:
\begin{enumerate}
\item Submodules of tame modules are tame.
\item  If the finitely generated group $G$ is an abelian extension of the abelian group $Q$,
then $G$ is finitely presented if and only if $A$ is tame as a $\Z[Q]$--module.
\end{enumerate}

\section{The geometry of the Cayley graph and metabelian presentations}

Gromov, in his paper in  his paper
\cite{Gromov:1983-1}, introduced into group theory a view of
infinite groups as geometric objects. What follows here is a very
short account of Gromov's basic idea and possible applications to
the study of finitely generated metabelian groups.

Suppose that $G$ is a group equipped with a finite set $X$ of
generators.  We turn $G$ into a metric space $\G=\G(G,X)$ by
defining the distance $d(g,h)$ between $g$ and $h$ to be the
length of the shortest $X$-word equal to $gh^{-1}$. It turns out
that $\G$ is an invariant of $G$, for if one changes the finite
generating sets involved, the metric spaces turn out to be
quasi-isometric (see \cite{Bridson-Haefliger:1999-1} for more details and
further explanation.) In the event that $\G$ has the geometry of a
two-dimensional hyperbolic space, then the group $G$ is said to be
hyperbolic. These hyperbolic groups can be described in terms of
the complexity of their word problems, i.e., in terms of what are
called isoperimetric functions. Roughly speaking,  given a finite
presentation of a group, the isoperimetric function counts the
number of times one uses defining relations to prove that a given
relator of the group is the identity.This function is often
referred to as the Dehn function. Hyperbolic groups turn out to be
those finitely presented groups with linear Dehn functions. Little
is known about isoperimetric functions of finitely presented
solvable groups. Recently Kassabov and Riley~\cite{Kassabov-Riley:2011-1} showed that the example described in $\mathsection 11$, below, has an exponential
Dehn function, but when the extra relation $b^n=1$ is added, the
Dehn function becomes quadratic. This fascinating area needs further investigation. It would be interesting
to explore whether any of Gromov's ideas can be adapted to the
study of finitely generated metabelian groups.

Because finitely generated metabelian groups are finitely presentable in the variety of metabelian groups one can associate with them a
corresponding relative isoperimetric function. The first named
author computed some of these isoperimetric functions of
metabelian groups in Rio many years ago. Cheng-Fen Fuh in her
thesis at the Graduate Center of Cuny in 2000, computed the
relative isoperimetric functions of a number of additional
metabelian groups. In particular she proved that the group
$\langle a,t\ \vert \ a^t=a^2\rangle$ has a linear relative
isoperimetric function. In the category of all groups, this group
has an exponential isoperimetric function. It is not clear how
much of Gromov's program can meaningfully be carried out for
metabelian groups, but it is conceivable that some of these ideas
can be used. This is a topic that has yet to be
explored.

\section{Hilbert functions}

Let $G$ be a finitely generated metabelian group,
$\gamma_n(G)$ the n-th term of the lower central series of $G$.
The Lazard Lie ring $L(G)$ of $G$ is defined to be the additive
abelian group
\[
L(G)=\bigoplus_{n=1}^{\infty} \gamma_n(G)/\gamma_{n+1}(G)
\]
with binary operation
\[
[a\gamma_{m+1}(G),b\gamma_{n+1}(G)]=[a,b]\gamma_{m+n+1}(G).
\]
$L(G)$ can be thought of as a linearization of $G$. Then it turns
out that L(G) is metabelian and satisfies the maximal condition on
ideals. Several authors have proven, following Hilbert and Serre,
that $h(G)=\Sigma_{n=1}^{\infty}r_nt^n$, where $r_n$ is the
torsion-free rank of $\gamma_n(G)/\gamma_{n+1}(G)$, is a rational
function, termed the Hilbert function. One can define other growth
functions for these metabelian groups following some of the ideas
that arise in  algebraic geometry.

\section{Where to go from here}

As already noted, the major open problem about
finitely generated metabelian groups is the isomorphism problem.
The structure of finitely generated metabelian groups  depends on
the ideal theory
 of finitely generated commutative rings, and the corresponding submodule theory of
the modules over these rings.  Adjunction of inverses to rings and modules, that is, localization, results in simpler rings and modules. Localization allows one to focus attention on rings
with exactly one maximal ideal, the so-called local rings.  These rings
need  not be finitely generated, but they remain Noetherian, as do
the modules attached to them. This powerful tool,  localization,  has been much used in algebraic geometry.

Localization, and the related tool of completion, provide
 the means for distinguishing one algebraic variety from another.
 They also  provide an important way, not an easy one,  to study finitely generated metabelian groups
and in particular a possible way to attack the isomorphism
problem.  The major object of this paper is to sketch some of the
ideas that will be used in our ongoing work to investigate this
problem.  We focus on the flexibility provided by localization to
help deepen our understanding of finitely generated metabelian
groups.  Note, for instance, that the first named author used
localization in the proof that finitely generated metabelian
groups can be embedded in finitely presented metabelian groups.
This will be discussed further in the sequel. However in order to
better understand the point of view that we have used in
investigating finitely generated metabelian groups, we first
recall some notions from algebraic geometry.

\section{Affine algebraic sets, localization and completions}

The model used in our study of finitely generated
metabelian groups is an affine algebraic set. We remind the reader
that an affine algebraic set is a set $X$ of points in $k^n$,
where $k$ is a field, consisting of the zeroes of a set of
polynomials in $k[x_1,\dots,x_n]$.  These algebraic sets define a
topology on $k^n$, in which they are the closed sets. This
topology is known as the Zariski topology.
 An irreducible algebraic set, $X$, is called a variety, and an algebraic set $X$ is a variety
if the set of polynomials that vanish on $X$ is a prime ideal.
The set $A=A(X)$ of polynomial functions on $X$ can be identified with
$k[x_1,\dots,x_n]/I(X)$, where $I(X)$ is the set of polynomials vanishing on $X$.
$A$ is called the coordinate ring of $X$.

Among the points $P$ on $X$ are the so-called non-singular ones. These are the points at which at
least one of the partial derivatives does not vanish at  $P$.
The set of functions in $A$ which vanish at such a non-singular point $P$ is an
ideal, $\mathfrak{m}$, of $A$ and  $A/\mathfrak{m}\cong k$.
It follows that $\mathfrak{m}$
is a maximal ideal of $A$ and so prime.
 $S=A-\mathfrak{m}$ is multiplicatively closed.  Adjoin to $A$ inverses of the elements of $S$ - we
 will describe this construction in more detail in  $\mathsection 10$.
The result is a ring denoted $A_S$
containing $A$ in which the elements of  $S$
are now invertible.
 Adjoining the inverses of the elements of $S$ to $\mathfrak{ m}$, denoted $\mathfrak{m}_S$,  yields a maximal ideal of
 $A_S$, the unique
maximal ideal of $A_S$.
 So $A_S$ is a local ring and is called the local ring of the point
$P$.
 $A_S$ is called the field of fractions of this local ring.
In a local ring $A$ with  maximal ideal $\mathfrak{m}$ the powers
of $\mathfrak{m}$ have intersection \{0\}. It follows that   $A$
embeds in its completion using the powers of $\mathfrak{m}$. Two
such local rings are isomorphic if and only if the varieties on which
the given non-singular points lie are birationally equivalent. The
latter holds only if the two fields of fractions are isomorphic.
These ideas can be adapted to study finitely generated metabelian
groups.

\section{A new approach to the isomorphism problem for finitely generated metabelian groups}

It seems that some of the ideas described above for
studying the birational equivalence of affine algebraic sets can
be carried over to studying finitely generated metabelian groups.
To do so, we modify the techniques described above, namely
localization and completion.  In order to do so, we observe that a
finitely generated metabelian group $\Gamma$ contains a normal
subgroup $G$ of finite index which is residually nilpotent, see~\cite{Lennox-Robinson:2004-1}, p.73.  So one
can reduce the isomorphism problem for finitely generated
metabelian groups to a problem involving finite extensions of
finitely generated residually nilpotent groups and the isomorphism
problem for finitely generated residually nilpotent metabelian
groups. In the event that a finitely generated metabelian group
$G$ is residually nilpotent, both techniques, localization and
completion, become available. In this section we will describe how
this comes about.

To this end, we need to introduce some additional notation. Let
$G$ be a group. Then we define $\g_1(G)=G$ and inductively for
$n>1$, $\g_n(G)= [G,\gamma_{n-1}(G)]$. Then the series
\[
G=\g_1(G)\geq \g_2(G) \geq \dots \geq \g_n(G) \geq \dots
\]
is termed the lower central series of $G$.
$G$ is termed residually nilpotent if the intersection of the terms of the lower central series is trivial.
Under these circumstances, we can avail ourselves of a completion, by analogy with the one
we described in the discussion of algebraic geometry, the pro-nilpotent completion, the inverse
limit, $\hat{G}$, of the quotients $G/\g_n(G)$. Some of the properties of $G$ carry over to $\hat{G}$.  In
particular, if $G$ is polycyclic, then we will prove that the finitely generated subgroups of $\hat{G}$
are also polycyclic. We will be more interested in our attempts to better understand the isomorphism
problem, to investigate the implications when two groups $G$ and $H$ have the same lower central quotients, i.e. if
\[
G/\g_n(G)\cong H/\g_n(H),\,  \text{for every}\ n\geq 1 .
\]
If $G$ and $H$ are residually nilpotent groups then we term $H$
para-$G$ if there exists a homomorphism $\phi$ from $G$ to $H$
which induces isomorphisms from $G/\g_n(G)$ to $H/\g_n(H)$ for $n> 0$. This
notion will play an important role in our work.  One of the objectives of our work
is to prove that there are surprisingly close connections between $H$ and $G$
given the right hypothesis, as we will show in Theorem 10.2, our Telescope theorem,
and Theorems 12.1, 12.2 and 12.3. As a simple sample of such a connection we
note first the easily proved
\begin{theorem} Suppose that $G$ and $H$ are residually nilpotent metabelian groups.
If $H$ is para-$G$ and if $H$ is finitely generated, then $G$ is also finitely generated.
\end{theorem}
 It seems likely
that given two residually nilpotent metabelian groups, one can
decide algorithmically whether they have the same lower central
quotients. Of course, if $H$ is para-$G$, then $G$ and $H$ have the
same lower central quotients. If $G$ is free in some variety and $G$
and $H$ have the same lower central quotients, then $H$ is para-$G$.
In particular, if a group $H$ has the same lower central quotients as
the free group $G$, then $H$ is para-$G$. These para-$G$ groups
are called parafree. There exists finitely generated parafree
groups which are not free.  In fact, it appears that for most groups $G$, there are para-$G$ groups $H$ which are not isomorphic to $G$~\cite{Baumslag:1967-1,Baumslag:1969-1}.  In order to obtain some deeper connections between a group $H$ which is para-$G$,
and $G$, it turns out that localization can play a key role. Indeed, the use of
localization in the study of finitely generated residually nilpotent metabelian groups
parallels that of localization that arises in algebraic geometry detailed above. The
basic idea here goes back to Levine although the work of Baumslag
and Stammbach seem related ~\cite{Levine:1989-2,Levine:1989-1,Baumslag-Stammbach:1977-1}.

We note in passing that a related construction to localization, which like localization is a smaller version of completion, arose from work in low dimensional topology.  Jerome Levine defined what he called the {\em algebraic closure of a group} to study concordance of knots~\cite{Levine:1989-1,Levine:1989-2}. Let $f: G\to H$ be a group homomorphism which induces an isomorphisms of lower central quotients. Then $f$ induces an isomorphism of the pro-nilpotent completions. These huge groups, often uncountable, can be unwieldy. Levine observed that a smaller, countable subgroup of the pro-nilpotent completion effectively replaces completion when studying finitely generated groups. Homomorphisms that induce isomorphisms on lower central series quotients become isomorphisms after taking closure, and the closure can be described as a direct limit through para-equivalencies of groups. The group closure functor $G \mapsto \overline G$ that appeared in Levine's work has played a significant role in low dimensional topology.

For metabelian groups, the localization tools which have their
origins in commutative algebra, and are suggested by algebraic
geometry, agree with Levine's topologically inspired group
closure.

\section{Localization of finitely generated, residually nilpotent,  metabelian groups and the Telescope Theorem}

The objective of this section is to describe, in part, the role localization plays in our efforts.  Statement (3) below is our {\em Telescope Theorem}, which shows that in a sense, our group localization of a group $G$ has the property of being {\em locally $G$}.  That is, every finitely generated subgroup is contained in an isomorphic copy of the given group $G$.  This has substantial consequences, and motivates our forthcoming definition of para-equivalence of groups.

We will need  what can be viewed as an instantiation of the change of coefficients
involving the second cohomology of a group, which we term here  {\em an extension pushout}.
In order to define such an
extension pushout, suppose that
$H$ is a group and $B$ is an abelian normal subgroup of $H$ with quotient $P$. So $H$ is
the middle of a short exact sequence
$$0 \longrightarrow B \buildrel{\beta}\over{\longrightarrow}
H\longrightarrow P \longrightarrow 1.$$ Here $\beta$ is the
inclusion of $B$ in $H$. $H$ acts on $B$ as a group of
automorphisms by conjugation, as does $P$, and so we can form the
semi-direct product  $H\ltimes B$.  We will denote the elements of
$H\ltimes B$ simply as $hb\ (b\in B, h\in H)$.

Now suppose that $S = 1 + \ker{\{\Z[P] \to \Z\}}$, and that $s \in S$.  Consider a $\Z[P]$-module homomorphism $\gamma \colon B \to B$, given by $b \mapsto b\cdot s$. This allows us to form a new semi-direct product $H \ltimes B$.  Then it follows that the subgroup $K$ of  $H\ltimes B$ generated by the elements
$b^{-1}\gamma(b)\  ( b \in B)$ is normal in $H\ltimes B$. We call $(H\rtimes B)/K$ the
 {\em extension pushout of  $B$ into $B$ through $H$ via $\gamma$}. Then, given the conditions above,
  the  following lemma holds.
 \begin{lemma}
 \begin{enumerate}
 \item The extension pushout $E$ of  $B$ into $B$ through $H$ via $\gamma$ is isomorphic to $H$.
 \item  The canonical homomorphsm $\phi$ of $H$ into $E$ mapping $h\in H$ to the coset $hK$ in $E$
 is a monomorphism if $H$ is residually nilpotent, and to a proper subgroup if $\gamma$ is not onto.
\end{enumerate}
\end{lemma}
We will provide details of the proof of Lemma 10.1 in \cite{Baumslag-Mikhailov-Orr:2012-1}.

We come now to the formulation of our Telescope Theorem,
Theorem 10.2. To this end, let $G$ be a finitely generated, residually nilpotent, metabelian group, and
let $A$ be the derived group of $G$, written
additively as usual,  and let $Q=G/A$. Let $I$ be the augmentation
ideal  of the group ring $R$ of $Q$, i.e., the elements of $R$
with coefficient sum zero.  $A$ can be viewed as an $R$-module as
already noted. It then follows from the residual nilpotence of
$G$, that
\[
\bigcap_{n=1}^{\infty}A I^n =0.
\]
Now set $S=1+I$. Then $S$ is multiplicatively closed and contains the identity $1$.
Observe that if $a \in A, a \neq 0$ and $s\in S$, then $as \neq 0$. Suppose the contrary. Now
$s=1-\alpha$, $\alpha \in I$, and $0=as=a(1-\alpha)$ implies $a=a\alpha$ and therefore
$a=a\alpha^n$ for every $n$. Hence $a=0$. We now form what is termed $A$ {\em localized at}  $S$, which
we denote here by $A_S$,  and which consists of the equivalence classes of $A \times S$
with respect to the equivalence relation $\sim$ defined as follows:
\[
(a,s)\sim(b,t) \text{ if } (at-bs)u=0 \text{ for some } u \in S.
\]
We denote the equivalence class of $(a,s)$ by $a/s$ and turn $A_S$ into an $R$-module in the obvious way.

Since  $S$ is countable, we can enumerate the elements of  $S$:
\[
S=\{s_0 = 1, s_1, s_2, \dots, s_n, \dots\}.
\]
Now for each $i = 0, 1, \dots$, define $A_i = A$, and define the $\Z[Q]$-module homomorphism $\gamma_n \colon A_n \to A_{n+1}$ given by $a \mapsto a\cdot s_n$.  The homomorphism
\[
\psi \colon A_0 \to A_S \text{ defined by } a \mapsto a/1
\]
extends to a homomorphism $A_1 \to A_S$, and this inductively defines injective homomorphisms $A_n \to A_S$. Then $A_S$ becomes a properly ascending union of its submodules $A_n$.

\begin{theorem} ({\bf The Telescope Theorem})  Let $G$ be an extension of $A$ by $Q$.  Then $G$ is the middle term in a short exact sequence of groups
$$
A \buildrel{\phi}\over\hookrightarrow G \twoheadrightarrow Q
$$
where $\phi$ is the inclusion of $A=A_0$  in $G$. Now put $G_0=G$.
Notice that $A_i$ is a proper submodule of $A_{i+1}$. Now let
$G_1$ be the extension pushout of $A_0$ through $A_1\rtimes G_0$
via $\gamma_0$, $G_2$  the extension pushout of  $A_1$ through
$A_2 \rtimes G_1$  via $\gamma_1$, $\dots$, $G_n$ the extension
pushout of $A_{n-1}$ through $A_n \rtimes G_{n-1}$ and so on. Then
by Lemma 10.1 each of the $G_n$ is isomorphic to $G$ and the image
of $G_{n-1}$ under the canonical homomorphism of $G_{n-1}$ into
$G_n$ is a proper subgroup of $G_n$. If we abuse the naming of the
various groups involved in the sequence of groups above we obtain
a properly ascending sequence of groups
\[
G=G_0 <  G_1< \dots  <  G_n <  \dots <   G_{\infty}=\bigcup_{n\geq 0}G_n.
\]
If we now define  $G_S$ to be the extension pushout of  $A=A_0$ through $A_S\rtimes G$ via the map
$\psi$ defined above.  Then it turns out that
$G_S$ is isomorphic to $G_{\infty}$. We refer to $G_S$ as {\em $G$ localized at $S$}, or as
\em{the telescope of $G$}.
\end{theorem}

We will give a complete proof of the Telescope Theorem in  \cite{Baumslag-Mikhailov-Orr:2012-1}.
As previously stated, this construction coincides with Levine's
group closure for metabelian groups. The construction of $G_S$ and
$A_S$ from $G$ and $A$ is functorial.

\section{Embedding groups and localization}

It turns out that localization not only provides a
new tool for studying finitely generated metabelian groups, as we
show in the sections that follow, but also allows for new
embedding theorems which we hope will provide insights into this
amazingly complex family of groups. Here is an illustration of one
use of localization which sheds light on how one can embed a
finitely generated metabelian group in a finitely presented
metabelian group. We use a different multiplicative set than used
in most of this paper, and to avoid possible confusion, we denote
this new multiplicative set by $K$, and preserve the notation $S$
in this paper for the multiplicative set $1 + \ker\{\mathbb Z[Q]
\to Z\}$ used in our group telescope.

To this end, consider the group
$$G =\langle b, s, t\ \vert \ [s, t] = 1, b^s = bb^t , [b, b^t ] = 1 \rangle,$$
It turns out that $G$ is actually the result of localization. To see how this comes about,
let $W$ be the wreath product of one infinite cyclic group by another,
$$W =\langle b \rangle \wr \langle t \rangle.$$
It is worth explaining what this means. Here $W = AQ$, where $Q
=\langle t \rangle$ and $A$ is the free abelian group, freely
generated by the conjugates of $b$ by the powers of $t$. $W$ is
not a finitely presented metabelian group. Indeed it can be
presented in the form
$$W=\langle b,t\ \vert \ [b,b^{t^i}]=1,\ i=1,2,\dots\rangle$$
but it cannot be finitely presented. The problem here is that
infinitely many of these commuting relations are needed in order
to present $W$. Now let $R$ be the integral group ring of $T$  and
let $K$ be the set of powers of $1 + t$. Then
$$R_K = \{r/s \ \vert \ r \in  R, s\in K\}$$
and
$$A_K = \{a/s \ \vert \  a \in A, s\in K\}.$$
$A_K$ is then an $R_K$-module and $U = gp(t,s = 1 + t)$ is a free
abelian group on $s$ and $1+t$ and it acts on $A_K$ by right
multiplication. So we can form the semi-direct product  $A_K
\rtimes U$ of $A_K$ and $U$. This then is our group $G$ .

Note how these finitely many relations suffice to make all of the conjugates
of $b$ by the powers of $t$ commute.  The following simple observations show how this
comes about. First we have $[b,b^t]=1$. So
$$1=[b,b^t]=[b,b^t]^s=[b^s,b^{ts}]=[b^s,b^{st}]=[bb^t,(b^s)^t]= [bb^t, b^tb^{t^2}]$$
and from this it follows that $[b,b^{t^2}]=1$, and similarly, that
all of the conjugates of $b$ by the powers of $t$ commute. It is
this kind of trick that enables one to embed finitely generated
metabelian groups in finitely presented metabelian groups
\cite{Baumslag:1973-1}, \cite{Remeslennikov:1973-1}.

If we had chosen $K$ to be the set of non-zero divisors of $R$, and
$U$ to be the multiplicative subgroup of $R_K$ generated by $K$,
then the semi-direct product $U \ltimes B_K$ is a metabelian group with
an extremely interesting array of subgroups. We will use variants
of this kind of group in constructing completions of finitely
generated, residually nilpotent, metabelian groups.

A careful examination of  the group $G_S$ reveals  surprising aspects which
resemble the nature of the group $G$ discussed above. In fact it  allows us to
construct a slightly different way
of embedding a finitely generated metabelian group into a finitely presented
metabelian group. Although this does not differ substantially from the first named author's
original proof,  it  does lead to a slightly different view of that theorem. Indeed using the
ideas involved in the discussion of theTelescope Theorem, Theorem 10.1, enables
us to prove the following
new theorem which is not easily proved alongs the lines of the earlier proof of the general embedding theorem.

\bigskip\noindent
{\em Every finitely generated residually nilpotent metabelian group can be embedded in
a finitely presented residually nilpotent metabelian group}.

\bigskip\noindent

So the use of localization might well lead to a deeper understanding
of finitely generated metabelian groups and their close connections with finitely
presented metabelian groups.

\section{Some of the implications of the use of localization}

Now suppose that $G$ and $H$ are finitely generated
residually nilpotent metabelian groups and that $H$ is para-$G$.
Then $H$ comes equipped with a homomorphism $\phi \colon G \to H$
which induces an isomorphism on the corresponding quotients of their
lower central series. The homomorphism $\phi$ is  a monomorphism
and in particular induces an isomorphism from $G_{ab}=G/[G,G]$ to
$H_{ab}=H/[H,H]$ . So we can identify  $H_{ab}$ with $G_{ab}$
which we denote simply by $Q$. Let $R$ be the integral group ring
of $Q$, $I$ the augmentation ideal of $R$ and let $S=1+I$, as
before. Then we have this crucial results linking the telescopes
$G_S$ and $H_S$ and their submodules and some consequences of them
as detailed below.
\begin{enumerate}
\item The monomorphism $\phi$ induces an isomorphism
\[
\phi_S:  G_S \longrightarrow H_S.
\]
\item If we now put $A=[G,G]$ and  $B=[H,H]$ then $\phi_S$ induces an isomorphism between the $R$-modules
$A_S$ and $B_S$ and so $\phi_S^{-1}$ maps $B_S$ to $A_S$. In particular $\phi_S^{-1}$ maps $B$ into $A_S$.
\item Now $B$ is a finitely generated $R$-module and $A_S$ is an ascending union of its
$R$-submodules $A_i$ (see $\mathsection 10$), so there exists an integer $\ell$ such that $\phi_S^{-1}$
maps $B$ into $A_{\ell}$. As already noted, the $R$-submodules $A_i$ are isomorphic to the $R$-module $A$.
If now $G$ is finitely presented, then $A$ and also the $A_i$
are tame $R$-modules and so too are their submodules (see $\mathsection 4$). So $B$ is a tame
$R$-module. Consequently $H$ is also finitely presented. Thus we have proved the rather
surprising theorem:

\begin{theorem} Suppose that $G$ and $H$ are finitely generated residually
nilpotent metabelian groups and that $H$ is para-$G$. Then $H$ is
finitely presented if and only if $G$ is finitely presented.
\end{theorem}

\item Recall once again that the Telescope Theorem states that a countable group $G_S$ is a union of an increasing, countable sequence of copies of groups isomorphic to $G$.  If a finitely generated group $H$ is para-$G$, then for any generator of $H$, there is some group $G_n$ in this filtration of $G_S$ that contains that generator of $H$.  Thus we have  unexpected  theorem:
\begin{theorem}
Suppose $G$ is a countable, residually nilpotent, metabelian group.  If a finitely generated group $H$ is para-$G$, then there are inclusions
$G \leq H$ and $H \leq G$, both inducing isomorphisms on lower central series quotients.
\end{theorem}

We call this {\em para-equivalence of groups}.  That is, $G$ and $H$ are para-equivalent if $H$ is para-$G$ and $G$ is para-$H$.
Restated, the above theorem states that if $H$ is finitely generated and para-$G$, then $G$ and $H$ are para-equivalent.

\item Observe, as well, that for similar reasons, the Telescope Theorem implies another theorem:

\begin{theorem}  Suppose that $H$ is a finitely generated, residually nilpotent, metabelian group such that $H$ is para-$G$.
If $G$ is polycyclic, then so too is $H$.  That is, a group para-equivalent to a polycyclic group is also polycyclic.
\end{theorem}

\item Finally we note that we can prove a variation of Theorem 12.3 which takes a somewhat different form.

\begin{theorem} If $G$ is a finitely generated subgroup of the pro-nilpotent completion of a residually
nilpotent metabelian polycyclic
group, then $G$ is polycyclic.
\end{theorem}

\end{enumerate}

\section{Classifying para-equivalent metabelian groups}

he following theorem is mostly a summary of some of
the results from the prior two sections, accumulating these
results for the classification theorem and examples that follow.
However, statement $i)$ below holds some new information.

\begin{theorem}\label{thm:summary}
Let $G$ and $H$ be residually nilpotent, metabelian groups with $H$ finitely generated.  Then
\begin{itemize}
\item[i)]  $H$ is para-equivalent to $G$ if and only if there is a homomorphism $G \to H$ inducing an isomorphism $G_S \to H_S$.
\item[ii)]  The function $G \mapsto G_S$ is functorial.  In particular, $G$ and $H$ are isomorphic if and only if there is a diagram as follows, where the vertical homomorphisms are inclusions:
\[
\begin{diagram}
\node{G} \arrow{e,t}{\cong}\arrow{s} \node{H}\arrow{s}\\
\node{G_S} \arrow{e,t}{\cong} \node{H_S}
\end{diagram}.
\]
\end{itemize}
In particular, the isomorphism class of the $\Z[G_{ab}]$-module $[G,G]_S$ is an invariant of the para-equivalence class of $G$.
\end{theorem}

The difference between statements $i)$ and $ii)$ create the foundation for our forthcoming Classification Theorem~\ref{thm:classify} for isomorphism classes of para-equivalent groups.

If $Ann([G,G])$ is the annihilator of $[G,G]$, viewed as a  module over the quotient ring $ \Z[G_{ab}]$, then
it becomes a faithful module over the ring $R=\Z[G_{ab}]/Ann([G,G])$, i.e., each non-zero element of $R$
acts non-trivially on $[G,G]$. Following our analogy with algebraic geometry, we call $R$
 the {\em coordinate ring of $G$.}  The coordinate ring, too, is an invariant of the para-equivalence class of $G$.

\begin{theorem}
Let $G$ be a finitely generated, residually nilpotent group, and suppose $H$ is a para-$G$ group. Then the
coordinate ring of $G$ is isomorphic to the coordinate ring of $H$.
\end{theorem}

We classify isomorphism classes of groups para-equivalent to $G$ using properties $i)$ and $ii)$ of Theorem~\ref{thm:summary}, and some basic observations about submodules of $[G,G]_S$.  Toward that goal, we first consider submodules $C \leq [G, G]_S$.  Following the theory of ideal classes in Dedekind domains, we call such a sub-module an {\em $S$-fractional submodule} if it is finitely generated and the inclusion of $C$ in $[G,G]_S$ induces an isomorphism
between $C_S$ and  $[G,G]_S$. We denote the set of $S$-fractional submodules by $\mathcal F([G,G]_S)$.

\begin{lemma}
Each $S$-fractional submodule $C \subset [G, G]_S$ determines a para-$G$ group.  We'll call this {\em the para-$G$ group determined by $C$.}
\end{lemma}

In order to see how to prove Lemma 13.3, we need to use our Telescope Theorem.  Since $[G, G]$ is a finitely generated module, there is some element $s \in S$ such that $s\cdot [G,G] \subset C$.  The resulting homomorphism $[G,G]\to C$
induces a homomorphism of second cohomology groups $H^2(G_{ab};[G,G]) \to H^2(G_{ab}; C)$ and thereby gives
rise to an extension of $C$ by $G_{ab}$.  The resultant extension is then what we have termed the para-$G$ group determined by $C$.
Note then an automorphism of $G_S$ determines an automorphism of $[G,G]_S$, and therefore an
action of $Aut(G_S)$ on $\mathcal F([G,G]_S)$. We term two fractional $S$-modules equivalent if such an induced
automorphism of $[G,G]_S$ maps one onto the other.

\begin{definition}
We define  the {\em ideal class monoid of $G$, $\mathcal C\ell(G)$} to be the set of equivalence classes
\[
\mathcal C\ell(G) = \frac{\mathcal F([G, G]_S)}{Aut(G_S)}.
\]
That is, $C\ell(G)$ is the set of $S$-fractional modules where two are equivalent if an induced automorphism of $[G,G]_S$ maps one onto the other.
\end{definition}

Some examples in the remaining sections, and especially $\mathsection 14$, will help to justify using the term {\em ideal class monoid.} Then we have the following  important theorem.

\begin{theorem}\label{thm:classify}
Isomorphism classes of groups para-equivalent to $G$ lie in one-to-one correspondence to elements of $\mathcal C\ell(G)$.
\end{theorem}

\section{Examples}

We give a number applications of
Theorem~\ref{thm:classify}, with minimal explanation.  Details
will appear in forthcoming papers.

For the remainder of this paper we always assume that $G$ is a finitely generated, residually nilpotent, metabelian group.
{\em These hypotheses will not be repeated.}

The coordinate ring of $G$ plays a central role in calculating the ideal class monoid $\mathcal C\ell(G)$.  In this section and the next, we consider the simplest cases, where the coordinate ring is a principal ideal domain or a Dedekind domain.  The module theory of Dedekind domains is highly structured, and should admit a deeper analysis than given here.

Most of our examples in this section lie in a rich, but on the face of it, a relatively simple class of groups.
Throughout the following discussion, we  denote the infinite
cyclic group with generator $t$ by $T$ and, as usual,  use multiplicative notation for $T$. So the integral group ring
$\Z[T]$ of $T$ is simply the ring of polynomials in  $t$ and $t^{-1}$, i.e., consists of  finite Laurent polynomials
in $t$.  We assume that the groups $G$ that we  consider here
 take the form $T \ltimes A$,  where $A$ is an  abelian group, and
 that $G_{ab}$ modulo the torsion subgroup of $G_{ab}$ is infinite
 cyclic.  We denote the coordinate ring of $G$ by $R$ and by $Aut_R(A_S)$ the group of $R$-homomorphisms
 of $A_S$. Now $Aut_R(A_S)$  acts on $\mathcal F(A_S)$ which again gives rise to an equivalence
 relation on $\mathcal F(A_S)$. We denote the set of equivalence classes under this equivalence relation
 by  $\mathcal C\ell_S(A)$. It then turns out, for groups of the form $T \ltimes A$ as above, that

\begin{lemma}\label{lem:split}
\[
\mathcal C\ell(G) = \mathcal C\ell_S(A).
\]
\end{lemma}

If  $A$ is the additive group of the  coordinate ring $R$ of $G$ then the group $U$ of units in the ring $R$ acts on $A$
by right multiplication and we can form the semi-direct product $U \ltimes A$. Many of our examples will either take this form or variations of it.

Now Lemma~\ref{lem:split} states that for groups which split over an infinite cyclic group, the action of $Aut(G_S)$ on the set of $S$-fractional modules
is equivalent to the action of the group of module automorphisms of $A_S$.
The simplest case is that in which  the coordinate ring $R$ is a principal ideal domain.  We have the following theorem:

\begin{theorem}\label{thm:pid}
Let $A$ be the additive group  of the ring $R = \Z[T]/J$, where $J$ is
an ideal of $\Z[T]$ and let  $G = T \ltimes A$. If  $R$ is a
principal ideal domain, then any finitely generated para-$G$ group
is isomorphic to $G$.
\end{theorem}

For instance, consider the groups $G = T \ltimes \Z[1/n] $, $n \neq 2$, where we view  $\Z[1/n]$ is the additive
subgroup of the subring of the rational numbers with denominators a power of $n$.  (The example $n = 2$ is not residually nilpotent.)
Here  $t^n$ acts on $\Z[1/n]$ by multiplication by $n$.  Then the coordinate ring of $G$ is $\Z[1/n]$, a principal ideal domain.
Hence any para-$G$ is isomorphic to $G$ by Theorem~\ref{thm:pid}.

Another interesting example is the Lamplighter group,  the wreath
product of the group $C_2=\langle a \ \vert \ a^2=1\rangle$ of order 2 and
the infinite cyclic group $T=\langle t\rangle$:
\[
L =C_2\wr T= \langle a, t\   \vert \ a^2, [a, t^{-k}at^{k}] \text{ for all $k$ } \in \Z\rangle.
\]
This coordinate ring of $L$ is $\Z[T]$,  a Laurent polynomial ring over a field, and thus a principal ideal domain.  Again Theorem~\ref{thm:pid} applies.

One can  construct numerous other interesting examples.   Here  we give an example of a group $G$ such that any group para-equivalent to $G$ is isomorphic to $G$ but for which  the coordinate ring for $G$ is not a principal ideal domain.

\begin{example}
Consider the group $G = C_2\ltimes A$, where now $A$ is the additive group of the integral group ring of $C_2$
and $C_2$ acts on $A$ by right multiplication.
  Then the coordinate ring of $G$ is $\Z[C_2]$. One can show that an ideal $J \leq \Z[C_2]$ is $S$-fractional
  if and only if $J \cap S \not = \emptyset$.  But in the ring $\Z[C_2]$, every $S$-fractional ideal is a principal ideal, and hence, every group para-equivalent to $G$ is isomorphic to $G$.
\end{example}

We conclude with an example where $\mathcal C\ell(G)$ is infinite.  We give only one example.
To this end,  consider the wreath product:
\[
W=\langle b\rangle  \wr \langle t\rangle = \langle t, b \ \vert \
[b, t^{-k}bt^{k}] \text{ for all $k$ } \in \Z \rangle.
\]
Then, for example, if $J$ is the ideal of $\Z[T]$ generated by  $ 2t-1$ and $2-t$.
The groups $ T \ltimes  J  $ and $W$ are para-equivalent. However
$J$ is not a principal ideal and therefore  the two groups are not
isomorphic. Contrast this with the Lamplighter example, a quotient
group of the above example, and one for which para-equivalence
implies isomorphism.

\section{Further examples and Dedekind domains}

We consider the next simplest example, where the
coordinate ring is a Dedekind domain.

Recall that a Dedekind domain is an integral domain such that non-zero proper ideals factor uniquely into products of prime ideals.  Such rings are integrally closed, i.e., they contain the roots of monic polynomials over the integers. The ring of algebraic numbers in a number field is a Dedekind domain, as well as the coordinate ring of a nonsingular, geometrically integral, affine algebraic curve over a field $\Bbbk.$

We will consider special Dedekind domains, that is, Dedekind domains which are quotient rings of a Laurent polynomial ring on a single variable.  We call such a Dedekind domain a {\em Laurent domain}.

We first consider examples of the form $G =T \ltimes D$ where $D$ is  a ring of algebraic integers over a real quadratic extension of the rationals.  More precisely we consider fields of the form $\Q(\sqrt{d})$ where $d$ is a square free positive integer.  In any such number field, the ring of algebraic integers is the field of the roots of monic polynomials with coefficients in $\Z$.  This is always a Dedekind domain, and for quadratic extensions has the form
\[
D = \Z[\sqrt{d}] \text{  for  } d \equiv_4 2, 3 \quad \text{ and } \quad D = \Z\left[\frac{1 + \sqrt{d}}{2}\right] \text{ for } d \equiv_4 1.
\]

Such rings have a norm which can be used to compute units in the ring.  Given a unit in the ring $D$, this defines a homomorphism $T \to D$ by sending $t$ to that unit. $D$ is Laurent if and only if $\iota$ is onto for some choice of unit. We made this computation for all $d < 100$, and the ring of algebraic integers in $\Q(\sqrt{d})$ is a Laurent domain for the following values of $d$:
\[
d = 2, 3, 10, 13, 15, 23, 26, 29, 35, 53, 77, 82, 85.
\]
The $S$-ideal class monoid $\mathcal C\ell_S(D)$ forms an abelian group when $D$ is a Dedekind domain, and elements in this group lie in one-to-one correspondence to isomorphism classes of groups para-equivalent to
$T \ltimes D$, by Lemma~\ref{lem:split}.  In these cases, the coordinate ring of  $T \ltimes D$
is precisely $D$.

The ring $D$ is a principal ideal domain for $d = 2, 3, 13, 23, 29, 53, \text{ and } 77$.

For the remaining cases above, one easily computes these abelian groups using elementary techniques from number theory, and there are exactly two isomorphism classes of groups in each para-equivalence class for each remaining value of $d$ above, with the exception of $d = 82$.  In this last case there are $4$ such groups, and the $S$-ideal class group is cyclic of order four.  For every Laurent domain we have computed, the homomorphism $\mathcal C\ell_S(D) \to \mathcal C\ell(D)$ is an isomorphism, but we see little reason to believe this holds in general.

To illustrate these results concretely, we give two isomorphism classes of groups which are para-equivalent.

If $R$ is a ring, then we denote the ideal of $R$ generated by the elements $a,b,\dots$ by $(a,b,\dots)$.
Consider the ideal $J$ of   $A=\Z[T]/(t^2 - 6t -1)$ generated
by the image of $(3, t-2)$.  Standard arguments show that this is a non-principal ideal.  (It's square is principal.)  The ring homomorphism $\Z[T] \to \Z[\sqrt{10}]$ defined by $t \mapsto 3 + \sqrt{10}$ induces an isomorphism $A \cong \Z[\sqrt{10}]$, and sends the ideal $J$ to $(3, 1+\sqrt{10})$.  One easily shows that both $J$ and $A$ are free abelian groups of rank two.  We can write thse as such, and the inclusion $J \subset A$ is the homomorphism
\[
J = \Z^2 \xrightarrow{\tiny \quad \begin{bmatrix} 3 & -2\\ 0 & \ \ 1 \end{bmatrix} \quad} \Z^2 = A.
\]
Here, $t$ acts on the domain and range, respectively, via the matrices:
\[
\begin{bmatrix}
2 & 3 \\ 3 & 4
\end{bmatrix}
\text{ and }
\begin{bmatrix}
0 & 1 \\ 1 & 6
\end{bmatrix}.
\]
This determines a para-equivalence of non-isomorphic groups
\[
T\ltimes J  \to T \ltimes A .
\]

We have the following theorem:
\begin{theorem}
For $D$ a Laurent domain, the $S$-ideal class group $\mathcal C\ell_S(D)$
is finite, and a subgroup of the usual ideal class group of $D$.
\end{theorem}

Another class of interesting Laurent domains arise from the ring of cyclotomic integers, $\Z[\zeta_n] \cong \Z[t, t^{-1}]/(\phi_n(t))$, where $\phi_n(t)$ is the n-th cyclotomic polynomial.  These rings are also Dedekind domains, in fact, the ring of algebraic integers in $\Q(\zeta_n)$, and hence, Laurent domains.
\begin{theorem}
Let $G =T \ltimes  \Z[\zeta_{n}] $, where the action of a generator $t$ of $T$ on $\Z[\zeta_n]$ is multiplication
by $\zeta_n$.
\begin{enumerate}
\item [i)] $G$ is residually nilpotent if and only if $n = p^k$ for some prime $p$ and positive integer $k$.
\item[ii)] In each of these cases, isomorphism classes of groups para-equivalent to $G$ lie in one-to-one correspondence with
\[
\mathcal C\ell_S(D) \cong \mathcal C\ell(D).
\]
Recall that the first is the group of $S$-fractional ideals of $D$ under the operation of multiplication of ideals, and modulo principal ideals. The latter group, $\mathcal C\ell(D)$, is the classical ideal class group of $D$, that is, one considers all ideals in $D$.
\end{enumerate}
\end{theorem}

$D = \Z[\zeta_n]$ is a principal ideal domain for $n < 23$, and thus any group para-equivalent to
$G = T \ltimes \Z[\zeta_{p^k}] $ is isomorphic to $G$ for prime powers $p^k < 23.$  Exactly three remaining groups, $D \rtimes T$, are determined by their lower central series quotients.  These are $p^k = 25, 27,$ and $32$.

The first interesting case occurs for $n = 23$.  In this case the following is a para-equivalence of non-isomorphic groups.
\[
T \ltimes \left(2,\frac{1 + \sqrt{-23}}{2}\right)  \subset T \ltimes \Z[\zeta_{23}].
\]

\bibliographystyle{amsalpha}

\bibliography{research}

\providecommand{\bysame}{\leavevmode\hbox to3em{\hrulefill}\thinspace}
\providecommand{\MR}{\relax\ifhmode\unskip\space\fi MR }
\providecommand{\MRhref}[2]{%
  \href{http://www.ams.org/mathscinet-getitem?mr=#1}{#2}
}
\providecommand{\href}[2]{#2}
\begin{thebibliography}{EFW07}

\bibitem[Bau67]{Baumslag:1967-1}
Gilbert Baumslag, \emph{Groups with the same lower central sequence as a
  relatively free group. {I}. {T}he groups}, Trans. Amer. Math. Soc.
  \textbf{129} (1967), 308--321. \MR{0217157 (36 \#248)}

\bibitem[Bau69]{Baumslag:1969-1}
\bysame, \emph{Groups with the same lower central sequence as a relatively free
  group. {II}. {P}roperties}, Trans. Amer. Math. Soc. \textbf{142} (1969),
  507--538. \MR{0245653 (39 \#6959)}

\bibitem[Bau73]{Baumslag:1973-1}
\bysame, \emph{Subgroups of finitely presented metabelian groups}, J. Austral.
  Math. Soc. \textbf{16} (1973), 98--110, Collection of articles dedicated to
  the memory of Hanna Neumann, I. \MR{0332999 (48 \#11324)}

\bibitem[BH99]{Bridson-Haefliger:1999-1}
Martin~R. Bridson and Andr{\'e} Haefliger, \emph{Metric spaces of non-positive
  curvature}, Grundlehren der Mathematischen Wissenschaften [Fundamental
  Principles of Mathematical Sciences], vol. 319, Springer-Verlag, Berlin,
  1999. \MR{1744486 (2000k:53038)}

\bibitem[BMOa]{Baumslag-Mikhailov-Orr:2012-1}
G.~Baumslag, R.~Mikhailov, and K.~E. Orr, \emph{Localization, completions and
  metabelian groups}.

\bibitem[BMOb]{Baumslag-Mikhailov-Orr:2012-2}
G.~Baumslag, R.~Mikhailov, and K.E. Orr, \emph{Ideal class theory and
  metabelian groups}.

\bibitem[BS77]{Baumslag-Stammbach:1977-1}
G.~Baumslag and U.~Stammbach, \emph{On the inverse limit of free nilpotent
  groups}, Comment. Math. Helv. \textbf{52} (1977), no.~2, 219--233.
  \MR{0463304 (57 \#3257)}

\bibitem[BS81]{Bieri-Strebel:1981-1}
Robert Bieri and Ralph Strebel, \emph{A geometric invariant for modules over an
  abelian group}, J. Reine Angew. Math. \textbf{322} (1981), 170--189.
  \MR{603031 (82f:20017)}

\bibitem[EFW]{Eskin-Fisher-Whyte:2007-2}
Alex Eskin, David Fisher, and Kevin Whyte, \emph{Coarse differentiation of
  quasi-isometries {II}: {R}igidity for {S}ol and {L}amplighter groups},
  no.~arXiv:0706.0940.

\bibitem[EFW07]{Eskin-Fisher-Whyte:2007-1}
\bysame, \emph{Quasi-isometries and rigidity of solvable groups}, Pure Appl.
  Math. Q. \textbf{3} (2007), no.~4, part 1, 927--947. \MR{2402598
  (2009b:20074)}

\bibitem[FM00]{Farb-Mosher:2000-1}
Benson Farb and Lee Mosher, \emph{Problems on the geometry of finitely
  generated solvable groups}, Crystallographic groups and their generalizations
  ({K}ortrijk, 1999), Contemp. Math., vol. 262, Amer. Math. Soc., Providence,
  RI, 2000, pp.~121--134. \MR{1796128 (2001j:20064)}

\bibitem[GM86]{Groves-Miller-Charles:1986-1}
J.~R.~J. Groves and Charles~F. Miller, III, \emph{Recognizing free metabelian
  groups}, Illinois J. Math. \textbf{30} (1986), no.~2, 246--254. \MR{840123
  (87i:20057)}

\bibitem[Gro84]{Gromov:1983-1}
Mikhael Gromov, \emph{Infinite groups as geometric objects}, Proceedings of the
  {I}nternational {C}ongress of {M}athematicians, {V}ol.\ 1, 2 ({W}arsaw, 1983)
  (Warsaw), PWN, 1984, pp.~385--392. \MR{804694 (87c:57033)}

\bibitem[Hal54]{Hall:1954-1}
P.~Hall, \emph{Finiteness conditions for soluble groups}, Proc. London Math.
  Soc. (3) \textbf{4} (1954), 419--436. \MR{0072873 (17,344c)}

\bibitem[KR]{Kassabov-Riley:2011-1}
Martin Kassabov and Tim Riley, \emph{The {D}ehn function of {B}aumslag's
  metabelian group}, no.~arXiv:1008.1966.

\bibitem[Lev89a]{Levine:1989-1}
Jerome~P. Levine, \emph{Link concordance and algebraic closure. {I}{I}},
  Invent. Math. \textbf{96} (1989), no.~3, 571--592. \MR{91g:57007}

\bibitem[Lev89b]{Levine:1989-2}
\bysame, \emph{Link concordance and algebraic closure of groups}, Comment.
  Math. Helv. \textbf{64} (1989), no.~2, 236--255. \MR{91a:57016}

\bibitem[LR04]{Lennox-Robinson:2004-1}
John~C. Lennox and Derek J.~S. Robinson, \emph{The theory of infinite soluble
  groups}, Oxford Mathematical Monographs, The Clarendon Press Oxford
  University Press, Oxford, 2004. \MR{2093872 (2006b:20047)}

\bibitem[LS77]{Lyndon-Schupp:1977-1}
Roger~C. Lyndon and Paul~E. Schupp, \emph{Combinatorial group theory},
  Springer-Verlag, Berlin, 1977, Ergebnisse der Mathematik und ihrer
  Grenzgebiete, Band 89. \MR{0577064 (58 \#28182)}

\bibitem[Rem73]{Remeslennikov:1973-1}
V.R. Remeslennikov, \emph{On finitely presented groups}, Proc. Fourth All-Union
  Symposium on the Theory of Groups, Novosibirsk (1973), 164--169.

\bibitem[Rip82]{Rips:1982-1}
E.~Rips, \emph{Generalized small cancellation theory and applications. {I}.
  {T}he word problem}, Israel J. Math. \textbf{41} (1982), no.~1-2, 1--146.
  \MR{657850 (83m:20047)}

\bibitem[Tar49a]{Tartakovski:1949-3}
V.~A. Tartakovski{\u\i}, \emph{Application of the sieve method to the solution
  of the word problem for certain types of groups}, Mat. Sbornik N.S.
  \textbf{25(67)} (1949), 251--274. \MR{0033815 (11,493b)}

\bibitem[Tar49b]{Tartakovski:1949-1}
\bysame, \emph{The sieve method in group theory}, Mat. Sbornik N.S.
  \textbf{25(67)} (1949), 3--50. \MR{0033814 (11,493a)}

\bibitem[Tar49c]{Tartakovski:1949-2}
\bysame, \emph{Solution of the word problem for groups with a {$k$}-reduced
  basis for {$k>6$}}, Izvestiya Akad. Nauk SSSR. Ser. Mat. \textbf{13} (1949),
  483--494. \MR{0033816 (11,493c)}

\end{thebibliography}
\end{document}